\date{July 30, 2001}
\title{A percolation formula}
\author{Oded Schramm}
\documentclass[12pt,naturalnames]{article}
\usepackage{amsmath}
\usepackage{amsthm}
\usepackage{amsfonts}
\usepackage{graphicx}
\newif\ifhyper\IfFileExists{hyperref.sty}{\hypertrue}{\hyperfalse}
\ifhyper\usepackage{hyperref}\fi

\newif\ifdraft
\draftfalse
\def\note#1/{\ifdraft {\bf [#1]}\fi}

\newtheorem{theorem}{Theorem}

\newtheorem{lemma}[theorem]{Lemma}

\def\eref#1{(\ref{#1})}
\def\QED{\qed\medskip}

\newcommand{\R}{\mathbb{R}}
\newcommand{\C}{\mathbb{C}}

\def\H{\mathbb{H}}
\def\U{\mathbb{U}}

\def\Im{{\rm Im}}
\def\Re{{\rm Re}}
\def\SLEkk#1/{$\mathrm{SLE}_{#1}$}
\def\SLEk/{\SLEkk{\kappa}/}
\def\SLE/{$\mathrm{SLE}$}
\def\Ito/{It\^o}

\def \P {{\bf P}}
\def\md{\mid}
\def\Bb#1#2{{\def\md{\bigm| }#1\bigl[#2\bigr]}}
\def\BB#1#2{{\def\md{\Bigm| }#1\Bigl[#2\Bigr]}}

\def\Pb{\Bb\P}

\def\PB{\BB\P}

\def \p {{\partial}}
\def \E {{\bf E}}

\def\closure{\overline}
\def\ev#1{{\mathcal{#1}}}
\def \proof {{ \medbreak \noindent {\bf Proof.} }}
\def\proofof#1{{ \medbreak \noindent {\bf Proof of #1.} }}

\def\iw{\hat w}

\def\BSmirnovPerc{Smirnov}
\def\BSmirnovPercII{Smirnov2}

\def\BhigherFunc{MR15:419i}

\def\BLSWi{math.PR/9911084}

\def\BSchSLE{MR1776084}

\def\CardyFormula{MR92m:82048}
\def\RSslecont{math.PR/0106036}
\def\Watts{MR98a:82059}
\def\LanglandsEtAl{MR94e:82056}
\def\GrimmettBook{Grimmett:book}
\def\KestenBook{MR84i:60145}
\begin{document}
\maketitle

\begin{abstract}
Let $A$ be an arc on the boundary
of the unit disk $\U$.
We prove an asymptotic formula
for the probability that there is a percolation cluster $K$
for critical site percolation on the triangular grid
in $\U$ which intersects $A$ and such that $0$ is surrounded
by $K\cup A$.
\end{abstract}

Motivated by questions raised by Langlands et al~\cite{\LanglandsEtAl}
and by M.~Aizenman,
J.~Cardy~\cite{\CardyFormula,math-ph/0103018}
derived a formula for the asymptotic probability for the existence of a
crossing of a rectangle by a critical percolation cluster.
Recently, S.~Smirnov~\cite{\BSmirnovPerc} proved Cardy's formula
and established the conformal invariance of critical
site percolation on the triangular grid.
The paper~\cite{\BLSWi} has a generalization of Cardy's formula.
Another percolation formula, which is still unproven, was derived
by G.~M.~T.~Watts~\cite{\Watts}.  The current paper
will state and prove yet another such formula.
A short discussion elaborating on the general context of these results
appears at the end of the paper.

\medskip

\begin{figure}
\centerline{\includegraphics*[height=2.3in]{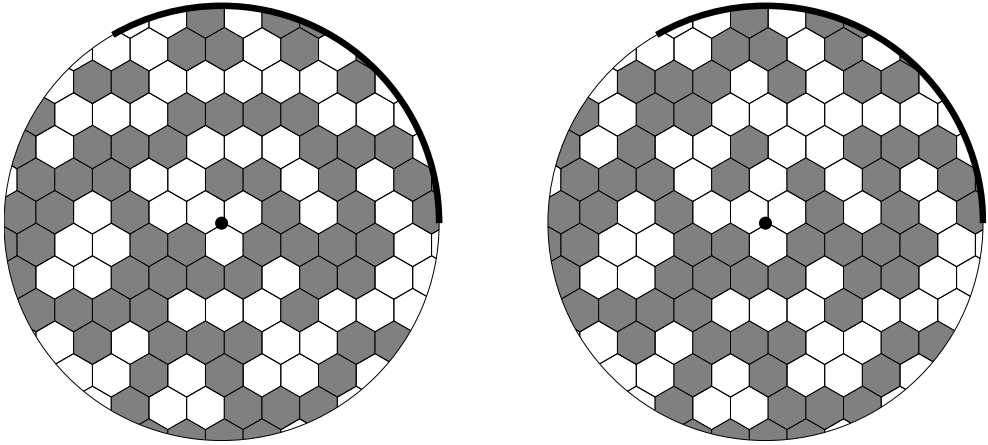}}
\caption{\label{f.twoways}The two topologically distinct ways to surround $0$.}
\end{figure}

Consider site percolation
on a triangular lattice in $\C$ with small mesh $\delta>0$,
where each site is declared open with probability $1/2$, independently.
(See~\cite{\GrimmettBook,\KestenBook} for background on percolation.)
It is convenient to represent a percolation configuration
by coloring the corresponding hexagonal faces of the dual grid;
black for an open site, white for a closed site.
Let $\mathfrak B$ denote the union of the black hexagons, intersected with
the closed unit disk $\closure{\U}$, and for
$\theta\in(0,2\pi)$ let $\ev A=\ev A(\theta)$
be the event that there is a connected component $K$
of $\mathfrak B$ which intersects the arc
$$
A_\theta:= \bigl\{e^{is}:s\in[0,\theta]\bigr\}\subset\p\U
$$
and such that $0$ is in a bounded component
of $\C\setminus\bigl(A_\theta\cup K\bigr)$.
Figure~\ref{f.twoways} shows the two distinct topological
ways in which this could happen.

\begin{theorem}\label{percform}
$$
\lim_{\delta\downarrow 0} \P[\ev A] =
\frac 12 -
\frac{\Gamma(2/3)}
{\sqrt{\pi}\,\Gamma(1/6)}
\,
  F_{2,1}\Bigl(\frac 12,\frac 23,\frac 32,-\cot^2\frac {\theta}{2}\Bigr)
\cot\frac{\theta}2
\,.
$$
\end{theorem}

Here, $F_{2,1}$ is the hypergeometric function.
See~\cite[Chap.~2]{\BhigherFunc} for background on hypergeometric functions.

\medskip
There is a second interpretation of the Theorem.
Suppose that $0$ is not on the boundary of any hexagon.
Let $C_1$ be the cluster of either black or white hexagons
which contains $0$.  Let $C_2$ be the cluster of the opposite
color which surrounds $C_1$ and is adjacent to it.
Inductively, let $C_{n+1}$ be the cluster surrounding and adjacent to $C_n$.
Let $m$ be the least integer such that $C_m$ is not contained
in $\U$, and let $C_m'$ be the component of
$\closure\U\cap C_m$ which surrounds $0$.  Let
$X:=1$ if $C_m'\cap\p\U\subset A_\theta$, let
$X:=0$ if $C_m'\cap A_\theta=\emptyset$, and otherwise
set $X:=1/2$.  Then $\P[\ev A]=\E[X]$.
This is so because
$$
\ev A = \{X=1\} \cup\{X=1/2\hbox{ and $C_m$ is black}\}\,.
$$

\medskip

Theorem~\ref{percform}
will be proved by utilizing the relation between the scaling
limit of percolation and Stochastic Loewner evolution with
parameter $\kappa=6$ (a.k.a.\ \SLEkk6/),
which was conjectured in~\cite{\BSchSLE} and proven by
S.~Smirnov~\cite{\BSmirnovPerc}.

We now very briefly review the definition and the relevant properties
of chordal \SLE/.
For a thorough treatment, see~\cite{\RSslecont}.
Let $\kappa\ge 0$, let $B(t)$ be Brownian motion on $\R$ starting
from $B(0)=0$, and set $W(t)=\sqrt{\kappa}\, B(t)$.
For $z$ in the upper half plane $\H$ consider the time flow $g_t(z)$ given by
\begin{equation}\label{chordal}
\p_t g_t(z)=\frac 2 {g_t(z)- W(t)}\,,\qquad g_0(z)=z\,.
\end{equation}
Then $g_t(z)$ is well defined up to the first time $\tau=\tau(z)$ such
that $\lim_{t\uparrow \tau} g_t(z)-W(t)=0$.
For all $t>0$, the map
$g_t$ is a conformal map from the domain $H_t:=\bigl\{z\in\H:\tau(z)>t\bigr\}$
onto $\H$.
The process $t\mapsto g_t$ is called Stochastic Loewner evolution
with parameter $\kappa$, or \SLEk/.

In~\cite{\RSslecont} it was proven that at least for $\kappa\ne 8$
a.s.\ there is a uniquely defined continuous path
$\gamma:[0,\infty)\to\closure{\H}$,
called the {\bf trace} of the \SLE/, such that for every
$t\ge 0$ the set $H_t$ is equal to the unbounded component
of $\H\setminus\gamma[0,t]$.
In fact, a.s.\
$$
\forall t\ge 0,\qquad\gamma(t)=\lim_{z\to W(t)} g_t^{-1}(z)\,,
$$
where $z$ tends to $W(t)$ from within $\H$.
Additionally, it was shown that $\gamma$ is a.s.\ transient,
namely $\lim_{t\to\infty}\bigl|\gamma(t)\bigr|=\infty$,
and that when $\kappa\in(0,8)$ for every $z_0\in\H$ we have
$\Pb{z_0\in\gamma[0,\infty)}=0$.

Fix some $z_0=x_0+i\,y_0\in\H$.
Then we may ask if $\gamma$ passes to the right or to the left of
$z_0$, topologically.  (Formally, this should be defined in terms of
winding numbers, as follows.  Let
$\beta_t$ be the path from $\gamma(t)$ to $0$ which follows
the arc $|\gamma(t)|\p\U$ clockwise from $\gamma(t)$ to $|\gamma(t)|$
and then takes the straight line segment in $\R$ to $0$.
Then $\gamma$ passes to the left of $z_0$ if the winding number of
$\gamma[0,t]\cup \beta_t$ around $z_0$ is $1$ for all large $t$.)
Theorem~\ref{percform} will be established by applying the following with $\kappa=6$:

\begin{theorem}\label{whichside}
Let $\kappa\in[0,8)$, and let $z_0=x_0+i\,y_0\in\H$.
Then the trace $\gamma$ of chordal \SLEk/ satisfies
$$
\Pb{\gamma\hbox{ passes to the left of }z_0}
=
\frac 12 +
\frac{\Gamma(4/\kappa)}
{\sqrt{\pi}\,\Gamma\bigl(\frac{8-\kappa}{2\kappa}\bigr)}\,
\frac{x_0}{y_0}
\,
  F_{2,1}\Bigl(\frac 12,\frac 4\kappa,\frac 32,-\frac {x_0^2}{y_0^2}\Bigr)
\,.
$$
\end{theorem}

When $\kappa=2,8/3,4$ and $8$ the right hand side simplifies to
$1 + \frac{x_0\,y_0}{\pi\,|z_0|^2}-\frac{\arg z_0}{\pi}$,
$\frac12+\frac{x_0}{|z_0|}$,
$1-\frac{\arg z_0}{\pi}$ and $\frac 12$, respectively.
\medskip

Let $x_t:=\Re\,g_t(z_0)-W(t)$, $y_t:=\Im\,g_t(z_0)$,
and $w_t:=x_t/y_t$.

\begin{lemma}\label{trans}
Almost surely, $\gamma$ is to the left [respectively, right]
of $z_0$ if $\lim_{t\uparrow\tau(z_0)} w_t=\infty$
[respectively, $-\infty$].
\end{lemma}

\proof
Suppose first that $\kappa\in[0,4]$.  In that case, a.s.\ $\gamma$ is a simple
path and $\tau(z_0)=\infty$, by~\cite{\RSslecont}.
Let $r>0$ be much larger than $|z_0|$, and let $\tau_r$ be the
first time $t$ such that $|\gamma(t)|=r$.
Let $D_+(r)\subset\H$ be the bounded domain whose boundary
consists of $[0,r]\cup\gamma[0,\tau_r]$ and
an arc on the circle $r\p\U$, and let
$D_-(r)\subset\H$ be the bounded domain whose boundary
consists of $[-r,0]\cup\gamma[0,\tau_r]$
and another arc on the circle $r\p\U$.
Given $\gamma$ we start a planar Brownian motion $B$ from
$z_0$.  With high probability, $B$ will hit $\gamma[0,\tau_r]\cup\R$
before exiting the disk $r\p\U$, provided $r$ is large.
This means that if $z_0\in D_+(r)$, then $B$ is likely
to hit $\gamma[0,\tau_r]$ from within $D_+(r)$, or hit
$[0,r]$.  By conformal invariance of harmonic measure,
this means that the harmonic measure in $\H$
of $[W(\tau_r),\infty)$ from $g_{\tau_r}(z_0)$ is close to $1$
if $z_0\in D_+(r)$ and close to zero if $z_0\in D_-(r)$.
Hence, the harmonic measure in $\H$ of $[0,\infty)$
from $x_{\tau_r}+i\,y_{\tau_r}$ is close
to $1$ if $z_0\in D_+(r)$ and close to zero if $z_0\in D_-(r)$.
Therefore, $w_{\tau_r}$ is close to $\pm\infty$ depending on
whether $z_0\in D_{\pm}(r)$.  This proves the lemma in the
case $\kappa\in[0,4]$.

For $\kappa\in(4,8)$, the analysis is
similar.  The difference is that a.s.\ $\gamma$ is
not a simple path, $\tau(z_0)<\infty$,
and $z_0$ is in a bounded component of
$\R\cup\gamma\bigl[0,\tau(z_0)\bigr]$ (see~\cite{\RSslecont}).
Clearly, $z_0$ is not in a bounded component
of $\R\cup\gamma[0,t]$ when $t<\tau(z_0)$.
Hence, at time $\tau(z_0)$ the path $\gamma$ closes a loop
around $z_0$.  The issue then is whether this is a clockwise
or counter-clockwise loop.  An argument as above shows that
this is determined by whether $w_t\to\pm\infty$
as $t\uparrow\tau(z_0)$.
\QED

\proofof{Theorem~\ref{whichside}}
Writing~\eref{chordal} in terms of the real and imaginary parts gives,
$$
dx_t = \frac{2\,x_t\,dt }{x_t^2+y_t^2}-dW(t)\,,
\qquad
dy_t = -\frac{2\,y_t\,dt}{x_t^2+y_t^2}\,.
$$
\Ito/'s formula then gives,
\begin{equation}\label{dw}
dw_t = -\frac {dW(t)}{y_t} +\frac{4 \,w_t\,dt}{x_t^2+y_t^2}\,.
\end{equation}
Make the time change
$$
u(t)=\int_0^t \frac{dt}{y_t^2}\,,
$$
and set
$$
\tilde W(t)=\int_0^t \frac{dW(t)}{y_t}\,.
$$
Note that $\tilde W/\sqrt{\kappa}$ is Brownian motion as a function
of $u$.
{}From~\eref{dw}, we now get
\begin{equation}\label{dwu}
dw = -d\tilde W + \frac{4 \,w\,du}{w^2+1}\,.
\end{equation}
We got rid of $x_t$ and $y_t$, and are left
with a single variable diffusion process $w(u)$.
(This is no mystery, but a simple consequence of scale
invariance.)  Given a starting point $\iw\in\R$ for the diffusion~\eref{dwu},
and given $a,b\in\R$ with $a<\iw<b$, we are interested in the probability
$h(\iw)=h_{a,b}(\iw)$ that $w$ will hit $b$ before hitting $a$.
Note that $h(w_u)$ is a local martingale.
Therefore, assuming for the moment that $h$ is smooth, by \Ito/'s formula,
$h$ satisfies
$$
\frac {\kappa}2\, h''(w)+ \frac{4\, w}{w^2+1}\,h'(w)=0\,,
\qquad h(a)=0\,,\qquad h(b)=1\,.
$$
By the maximum principle, these equations have a unique solution, and
therefore we find that
\begin{equation}\label{his}
h(w)= \frac{f(w)-f(a)}{f(b)-f(a)}\,,
\end{equation}
where
$$
f(w):=  F_{2,1}(1/2,4/\kappa,3/2,-w^2)\,w\,.
$$
We may now dispose of the assumption that $h$ is smooth, because
\Ito/'s formula implies that the right hand side in~\eref{his}
is a martingale, and it easily follows that it must equal $h$.
By~\cite[2.10.(3)]{\BhigherFunc} and our assumption $\kappa<8$ it follows that
\begin{equation}\label{finitelim}
\lim_{w\to\pm\infty} f(w)
=
\pm
\frac{ \sqrt{\pi}\, \Gamma\bigl((8-\kappa)/(2\kappa)\bigr)}{2\,\Gamma(4/\kappa)}
\,.
\end{equation}
In particular, the limit is finite, which
shows that $\lim_{b\to\infty} h_{a,b}(w)>0$ for all $w>a$.
Hence, the diffusion process~\eref{dwu} is transient.
Moreover,
$$
\PB{\lim_{u\to\infty} w_u = +\infty}
= \frac{f(\iw)-f(-\infty)}{f(\infty)-f(-\infty)}\,.
$$
An appeal to the lemma now completes the proof.
\QED

\proofof{Theorem~\ref{percform}}
For simplicity, assume that the points $0$, $1$ and $e^{i\theta}$ are
not on the boundary of any hexagon in the grid.
As above, let $\mathfrak B$ be the intersection of
the union of the black hexagons with
$\closure{\U}$, and let $\hat{\mathfrak B}$ be the union of
$\mathfrak B$ and the set
$S:=\bigl\{r\,e^{is}:r\ge 1,\,s\in[0,\theta]\bigr\}$.
Let $\beta$ be the boundary of the intersection of $\closure{\U}$ with the
component of $\hat{\mathfrak B}$ containing $S$.
Then $\beta$ is a path in $\closure{\U}$ from $1$ to $e^{i\theta}$.
It is immediate that the event $\ev A$ is equivalent to
the event that $0$ appears to the right of the path $\beta$; that is,
that the winding number of the concatenation of $\beta$ with the
arc $A_\theta$ with the clockwise orientation around $0$ is $1$.

S.~Smirnov~\cite{\BSmirnovPerc}
has shown that as $\delta\downarrow 0$
the law of $\beta$ tends weakly to the law of the image of
the chordal \SLEkk6/ trace $\gamma$ under any fixed conformal map
$\phi:\H\to\U$ satisfying $\phi(0)=1$ and $\phi(\infty)=e^{i\theta}$.
(See also~\cite{\BSmirnovPercII}.)
We may take
$$
\phi(z)=
e^{i\theta}\,
\frac{z+\cot\frac {\theta}2-i}
{z+\cot\frac {\theta}2+i}\,.
$$
The theorem now follows by setting $\kappa=6$ in Theorem~\ref{whichside}.
\QED

\bigskip {\noindent\bf Discussion.} 
According to J.~Cardy (private communication, 2001),
presently, the conformal field theory methods used by him to
derive his formula do not seem to supply even a heuristic
derivation of Theorem~\ref{percform}.
On the other hand, it seems that, in principle, probabilities for 
``reasonable'' events
involving critical percolation can be expressed as solutions of
boundary-value PDE problems, via SLE$_6$.
But this is not always easy. 
In particular, it would be nice to obtain a proof of
Watts' formula~\cite{\Watts}.
The event $\ev A$
studied here was chosen because the corresponding proof is particularly
simple, and because the PDE can be solved explicitly.

\bigskip {\noindent\bf Acknowledgements.}
I am grateful to Itai Benjamini and to Wendelin Werner for
useful comments on an earlier version of this manuscript.


\end{document}